

\baselineskip=14pt
\parskip=10pt

\magnification=\magstephalf

\def\1{{\overline{1}}}
\def\2{{\overline{2}}}
\parindent=0pt
\overfullrule=0in

\def\frac#1#2{{#1 \over #2}}
\centerline
{\bf 
Efficient Evaluations of Weighted Sums over the Boolean Lattice 
}
\centerline
{\bf 
inspired by  conjectures of  Berti, Corsi, Maspero, and  Ventura
}
\bigskip
\centerline
{\it Shalosh B. EKHAD and Doron ZEILBERGER}
\bigskip

{\bf Abstract:} In their study of water waves, Massimiliano Berti, Livia Corsi, Alberto Maspero, and Paulo Ventura,
came up with two intriguing conjectured identities involving certain weighted sums over the Boolean lattice. 
They were able to prove the first one, while the second is still open.
In this methodological note, we will describe how to generate many terms of these types of weighted sums, and if in luck,
evaluate them in closed-form. We were able to use this approach to give a new proof of their first
conjecture, and while we failed to prove the second conjecture, we give overwhelming evidence for its veracity.
In this second version, we are happy to announce that Mark van Hoeij was able to complete the proof of the second conjecture,
by explicitly solving the second-order recurrence  mentioned at the end.

{\bf An  Intriguing Email message from Alberto Maspero}

Awhile ago one of us (DZ) received an email message [M],
with the following. 

\quad \quad {\it In our current study of water waves [BCMV], continuing our work in [BMV], and other papers, we came across the following sums.}

Let $p \geq 2$, $1 \leq q \leq p-1$ and $0<j_1 <j_2 < \dots <j_q <p$ be positive integers. Define
$$
n_q^{(p)}(j_1, \dots, j_q) :=
(-12)^q j_1 \cdots j_q \,
\frac{ j_1(j_2-j_1) \cdots (j_q-j_{q-1})(p-j_q)}
{(p^3-p j+ j_1-j_1^3) \cdots (p^3-p+j_q-j_q^3)} \quad,
\eqno(1)
$$
and
$$
g_q^{(p)}(j_1, \dots, j_q) :=
-\frac{28}{9}p^2 + \frac{49}{45} q + \frac{32}{9} j_1^2 -\frac{4}{9} -\frac{4}{9} \frac{p^3-p}{j_1}+ \frac{5}{18} (p^3-p) \left (\frac{1}{j_1} + \cdots + \frac{1}{j_q} \right )
\eqno(2)
$$
$$
\frac{38}{15} \left ( j_1^2 + \cdots + j_q^2 \right ) +
\frac{p^3-p}{5}  \left( \frac{p+j_1}{p^2+j_1^2+pj_1-1} + \dots + \frac{p+j_q}{p^2+j_q^2+pj_q-1} \right)
-\frac{13}{9} (j_1 j_2 + j_2 j_3 +\cdots + j_q p) \quad.
$$

The following identities seem to hold.
$$
\sum_{q=1}^{p-1} \,\,\, \sum_{0<j_1 < \dots <j_q <p} n_q^{(p)} (j_1, \dots , j_q) \, = \, -p \quad.
\eqno(3)
$$
$$
\sum_{q=1}^{p-1} \,\, \,\sum_{0<j_1 < \dots <j_q <p} n_q^{(p)} (j_1, \dots , j_q) \,  g_q^{(p)} (j_1, \dots. j_q) \, 
= \, p(p+1)^2 \quad.
\eqno(4)
$$

Maspero concluded:

{\it ``We actually managed, after dire efforts, to prove $(3)$, whereas claim $(4)$ seems out of reach. We just verified it for $p=2, \dots, 21$.''}

\vfill\eject

{\bf Why are such sums interesting?}

Note that the sums in $(3)$ and $(4)$ are {\it weighted sums} defined over all non-empty subsets of the set $\{1 ... p-1\}$. 
Hence to compute  many terms, straight from the definition, is exponentially expensive.

They are a bit unnatural in that
the symbol $p$ has both the roles of {\it integer} and of  {\it variable}.  Our first step is to {\it decouple} these two roles of $p$, introduce a {\it formal} variable
$x$, and define
$$
N(x;j_1, \dots, j_q) :=
(-12)^q j_1 \cdots j_q \,
\frac{ j_1(j_2-j_1) \cdots (j_q-j_{q-1})(x-j_q)}
{(x^3-x j+ j_1-j_1^3) \cdots (x^3-x+j_q-j_q^3)} \quad,
\eqno(1')
$$
$$
G(x;j_1, \dots, j_q) :=
-\frac{28}{9}x^2 + \frac{49}{45} q + \frac{32}{9} j_1^2 -\frac{4}{9} -\frac{4}{9} \frac{x^3-x}{j_1}+ \frac{5}{18} (x^3-x) \left (\frac{1}{j_1} + \cdots + \frac{1}{j_q} \right )
\eqno(2')
$$
$$
\frac{38}{15} \left ( j_1^2 + \cdots + j_q^2 \right ) +
\frac{x^3-x}{5}  \left( \frac{x+j_1}{x^2+j_1^2+xj_1-1} + \dots + \frac{x+j_q}{x^2+j_q^2+xj_q-1} \right)
-\frac{13}{9} (j_1 j_2 + j_2 j_3 +\cdots + j_{q-1} j_{q}+ j_q x) \quad.
$$

Note that $N(p;j_1, \dots, j_q)=n_q^{(p)}(j_1, \dots, j_q)$ and $G(p;j_1, \dots, j_q)=g_q^{(p)}(j_1, \dots, j_q)$.

We are interested in efficient computation, and if possible, explicit evaluation of
$$
A_p(x):=
\sum_{q=1}^{p-1} \sum_{0<j_1 < \dots <j_q <p} N(x;j_1, \dots , j_q) \, \quad,
\eqno(3')
$$
and
$$
C_p(x):=
\sum_{q=1}^{p-1} \,\,\, \sum_{0<j_1 < \dots <j_q <p} N (x;j_1, \dots. j_q)  G  (x;j_1, \dots. j_q) \, \quad .
\eqno(4')
$$

In order to facilitate {\it dynamical programming}, it is natural to consider these weighted sums where the largest member of the subset, $j_q$, is fixed. So we define

$$
B_p(x):=
\sum_{0<j_1 < \dots <j_q=p} \, N(x;j_1, \dots , j_q) \, \quad,
\eqno(3'')
$$

and
$$
D_p(x):=
\sum_{0<j_1 < \dots <j_q=p} \, N(x;j_1, \dots , j_q)  G(x;j_1, \dots , j_q) \, \quad .
\eqno(4'')
$$
Once the quantities $B_p(x)$ and $D_p(x)$ are known, our original quantities of interest, $A_p(x)$ and $C_p(x)$, can be evaluated using
$$
A_p(x)= \sum_{p'=1}^{p-1} B_{p'}(x) \quad.
$$
$$
C_p(x)= \sum_{p'=1}^{p-1} D_{p'}(x) \quad.
$$

{\bf The general framework}

Note that the weights $N(x; j_1, \dots, j_q)$ and  $G(x; j_1, \dots, j_q)$ have a {\it recursive ``Markovian'' structure}. 

$\bullet$ If you know  $N(x;j_1, \dots, j_{q-1})$, you can quickly get $N(x;j_1, \dots, j_{q-1}, j_q)$, by {\bf multiplying} by a certain function of $(j_{q-1},j_q)$.

$\bullet$ If you know  $G(x;j_1, \dots, j_{q-1})$, you can quickly get $G(x;j_1, \dots, j_{q-1}, j_q)$, by {\bf adding}  a certain (different) function of $(j_{q-1},j_q)$.

This leads us to consider the following general set-up.

{\bf Definition}: Let $f_1(X)$ be an arbitrary uni-variate function, and $f_2(X,Y)$ an arbitrary bivariate function. Define the weight, for singleton sets $\{j_1\}$
$$
W(f_1,f_2; [j_1]):=f_1(j_1) \quad,
$$
and for sets with more than one element, recursively (where we write $j_1, \dots, j_q$ in increasing order):
$$
W(f_1,f_2; [j_1, \dots, j_q]):= W(f_1,f_2; [j_1, \dots, j_{q-1}]) \cdot f_2(j_{q-1},j_q) \quad .
$$
Similarly let $g_1(X)$ and $g_2(X,Y)$ be arbitrary univariate and bivariate functions and define

$$
V(g_1,g_2; [j_1]):=g_1(j_1) \quad,
$$
and for sets with more than one element, recursively
$$
V(g_1,g_2; [j_1, \dots, j_q]):= V(g_1,g_2; [j_1, \dots, j_{q-1}])+ g_2(j_{q-1},j_q) \quad .
$$

Note that the original [BCMV] summations have the following $f_1,f_2,g_1,g_2$:

$$
f_1(X) \, = \, -\frac{12 X^{2} \left(x -X \right)}{-X^{3}+x^{3}+X -x} \quad,
$$
$$
f_2(X,Y) \, = \,
-\frac{12 Y \left(Y -X \right) \left(x -Y \right)}{\left(x -X \right) \left(-Y^{3}+x^{3}+Y -x \right)} \quad,
$$
$$
g_1(X)=
-\frac{28 x^{2}}{9}+\frac{29}{45}+\frac{274 X^{2}}{45}-\frac{x^{3}-x}{6 X}+\frac{\left(x^{3}-x \right) \left(x +X \right)}{5 X^{2}+5 x X +5 x^{2}-5}-\frac{13 x X}{9} \quad ,
$$
$$
g_2(X,Y)=\frac{\frac{5}{18} x^{3}-\frac{5}{18} x}{Y}+\frac{38 Y^{2}}{15}+\frac{\left(x^{3}-x \right) \left(x +Y \right)}{5 Y^{2}+5 x Y +5 x^{2}-5}-\frac{13 X Y}{9}-\frac{13 \left(Y -X \right) x}{9}+\frac{49}{45} \quad .
$$

The general analogs of  $A_p(x)$, $B_p(x)$, $C_p(x)$ and $D_p(x)$, let's call them $a_p$, $b_p$, $c_p$, and $d_p$, respectively: are
$$
b_p:=\sum_{0<j_1 < \dots <j_q=p} W(f_1,f_2; [j_1, \dots, j_q]) \quad,
$$
and then
$$
a_p = \sum_{p'=1}^{p-1} b_{p'} \quad.
$$
$$
d_p:=\sum_{0<j_1 < \dots <j_q=p} W(f_1,f_2; [j_1, \dots, j_q])\, V(f_1,f_2; [j_1, \dots, j_q]) 
\quad,
$$
and then
$$
c_p = \sum_{p'=1}^{p-1} d_{p'} \quad.
$$

Let's first try to examine $b_p$.

We can break-up the sum that defines $b_p$, where every summand has $j_q=p$, according to the value of $j_{q-1}$:
$$
b_p:=\sum_{{0<j_1 < \dots <j_q}  \atop {j_q=p}} W(f_1,f_2; [j_1, \dots, j_q]) \, = \,
\sum_{p'=1}^{p-1} \,\, \sum_{{0<j_1 < \dots <j_{q-1} <p}  \atop {j_{q-1}=p'}} W(f_1,f_2; [j_1, \dots, j_{q-1}, p]) \, = \,
$$
$$
\sum_{p'=1}^{p-1} \,\, \sum_{{0<j_1 < \dots <j_{q-1} }  \atop {j_{q-1}=p'}} W(f_1,f_2; [j_1, \dots, j_{q-1}]) \cdot f_2(p',p) \, = \,
$$
$$
\sum_{p'=1}^{p-1}   f_2(p',p) \left (\sum_{0<j_1 < \dots<j_{q-1}=p'} W(f_1,f_2; [j_1, \dots, j_{q-1}])\right) \, = \,
\sum_{p'=1}^{p-1}   f_2(p',p) \, b_{p'} \quad .
$$

Hence the sequence $b_p$ can be computed in quadratic-time using the recurrence
$$
b_p \, = \, \sum_{p'=1}^{p-1} f_2(p',p)\, b_{p'} \quad ,
$$
subject to the {\it initial condition}
$$
b_1=f_1(1) \quad .
$$

A similar argument, that we omit, enables us to get a quadratic-time  recurrence for $d_p$, that assumes that $b_p$  is already known.

Once we get a hold of $b_p$ and $d_p$, we can recover $a_p$ and $c_p$ using $a_p=\sum_{p'=1}^{p-1} b_{p'}$ and $c_p=\sum_{p'=1}^{p-1} d_{p'}$.

Let's go back to specializing to the [BCMV] (already proved by them) conjecture $(3)$.

If we are lucky and we can conjecture an explicit expression for $B_p(x)$, then all we have to do is verify that this conjectured expression also satisfies the same
recurrence and initial condition. With the above $f_2(X,Y)$ the recurrence becomes
$$
B_p(x) \, = \, -12\frac{p(x-p)}{x^3-x+p-p^3} \left ( p + \sum_{p'=1}^{p-1} \frac{p-p'}{x-p'} B_{p'}(x)  \right ) \quad ,
\eqno(5)
$$
with the initial condition $B_0(x)=0$.

Cranking out the first $20$ terms one easily conjectures
$$
B_p(x) \, = \, \frac{12 p^{2} \left(p -x \right)}{x \left(x +1\right) \left(x -1\right)} \quad,
$$
and it is routine to verify  (even by hand, but Maple is glad to do it for you) that $(5)$ is satisfied if $B_p(x)$ is replaced by the above right side.
Then we ask Maple to kindly sum
$$
A_p(x)= \sum_{p'=1}^{p-1}  \frac{12 p'^{2} \left(p' -x \right)}{x \left(x +1\right) \left(x -1\right)} \quad ,
$$
giving
$$
A_p(x) = {\frac {p \left( p-1 \right)  \left( 3\,{p}^{2}-4\,xp-3\,p+2\,x \right) }{x \left( x-1 \right)  \left( x+1 \right) }} \quad .
$$
Now what [BCMV] are really interested in  is not $A_p(x)$, in general,  but the special case $x=p$, i.e. in $A_p(p)$. Plugging-in $x=p$ above, and simplifying gives that indeed
$$
A_p(p) \, = \, -p \quad .
$$
So we have a new proof of the already-proved-by-them  identity $(3)$ of [BCMV].

We can get similar  dynamical programming (quadratic-time) recurrence for $D_p(x)$ that expresses it in terms  of  previous values  $\{D_{p'}(x): 1 \leq p' \leq p-1\}$ {\bf and} (the already known) $B_p(x)$. Using the above {\it proved}
expression for the latter, we can compute many terms. Alas, it is no longer a nice rational function, and the sequence seems very complicated. But using the
{\it holonomic ansatz} [Z] (see [K] for a great Mathematica implementation) one can first {\it guess} (very complicated!) linear recurrences for both $D_p(x)$ and $C_p(x)$,
that nevertheless, to our pleasant surprise, are mere second order (but with very complicated coefficients).
See procedures {\tt DxH(p,x)} and {\tt CxH(p,x)} in our Maple package mentioned below. 
These recurrences are first guessed, using {\it undetermined coefficients}, implemented in our Maple package {\tt FindRec.txt}, available from:

{\tt https://sites.math.rutgers.edu/\~{}zeilberg/tokhniot/FindRec.txt} \quad .

Once guessed, they are all automatically and {\it rigorously} provable using the holonomic ansatz as implemented by Koutschan, i.e. the sequences defined by these second-order recurrences
also satisfy the original recurrences.

This enables us to easily compute the first $2000$ terms of the sequence $\{C_p(x)\}$ that are all very complicated rational functions of $x$. But when
we plug-in $x=p$ the sequence $\{C_p(p)\}$ {\bf coincides} with the conjectured sequence $\{p(p+1)^2\}$. Of course, this is not
a rigorous proof, but being {\it empiricists},  knowing that it is true for the first $2000$ terms, is {\bf good enough for us}.

{\bf Maple package and input and output files}

Everything is implemented in the Maple package {\tt BCMV.txt} available from

{\tt https://sites.math.rutgers.edu/\~{}zeilberg/tokhniot/BCMV.txt} \quad.

The front of this article

{\tt https://sites.math.rutgers.edu/\~{}zeilberg/mamarim/mamarimhtml/bcmv.html} \quad,

contains the input and output files that {\it rigorously} prove the above explicit expressions for $A_p(x)$ and $B_p(x)$, and
that empirically verifies $(4)$ all the way to $p=2000$.

{\bf Postscript written Feb. 28, 2024}: Mark van Hoeij met our challenge, to explicitly solve the recurrence satisfied by $C_p(x)$, that enables
plugging-in $x=p$ into it and proving that indeed $C_p(p)=p(p+1)^2$.

See:

{\tt https://sites.math.rutgers.edu/\~{}zeilberg/mamarim/mamarimhtml/bcmvChallenge.txt} \quad.

For a detailed explanation, see the postscript kindly written by Mark van Hoeij:

{\tt https://sites.math.rutgers.edu/\~{}zeilberg/mamarim/mamarimhtml/bcmvMvH.html} \quad.

This completes the (rigorous!) proof of Conjecture $(4)$. A donation of $100$ dollars to the OEIS, in Mark van Hoeij's honor,  has been made.

{\bf References}

[BCMV] M. Berti, L. Corsi, A. Maspero, and P. Ventura, {\it Study of water waves}, in preparation.

[BMV] M. Berti, A. Maspero, and P. Ventura, {\it Full description of Benjamin-Feir instability of stokes waves in deep water},
 Inventiones mathematicae {\bf 230} (2022), 651-711. \hfill\break
{\tt https://link.springer.com/content/pdf/10.1007/s00222-022-01130-z.pdf} \quad , \hfill\break
arxiv version: {\tt https://arxiv.org/abs/2204.00809} \quad .

[K] Christoph Koutschan, {\it Advanced applications of the holonomic systems approach}, PhD thesis, \hfill\break
{\tt http://www.koutschan.de/publ/Koutschan09/thesisKoutschan.pdf} \quad .

[M] Alberto Maspero, {\it Email message to Doron Zeilberger}, Jan. 26, 2024.

[Z] Doron Zeilberger, {\it A Holonomic Systems approach to special functions}, J. Computational and Applied Math {\bf 32} (1990), 321-368. \hfill\break
{\tt https://sites.math.rutgers.edu/\~{}zeilberg/mamarim/mamarimhtml/holonomic.html} \quad .

\bigskip
\hrule
\bigskip

Shalosh B. Ekhad, c/o D. Zeilberger, Department of Mathematics, Rutgers University (New Brunswick), Hill Center-Busch Campus, 110 Frelinghuysen
Rd., Piscataway, NJ 08854-8019, USA. \hfill\break
Email: {\tt ShaloshBEkhad at gmail  dot com}   \quad .
\smallskip

Doron Zeilberger, Department of Mathematics, Rutgers University (New Brunswick), Hill Center-Busch Campus, 110 Frelinghuysen
Rd., Piscataway, NJ 08854-8019, USA. \hfill\break
Email: {\tt DoronZeil at gmail  dot com}   \quad .
\bigskip
{\bf First version: Feb. 22, 2024}   ; {\bf This second version: Feb. 28, 2024} . 

{\bf Exclusively published in the Personal Journal of Shalosh B. Ekhad and Doron Zeilberger and arxiv.org}

\end